\documentclass{article}

\usepackage{amsmath,amssymb,amsthm}

\usepackage{xspace}

\usepackage{setspace}

\newcommand{\defin}[1]{{\em #1}}
\newcommand{\C}{\ensuremath{{\mathbb C}}\xspace}
\newcommand{\R}{\ensuremath{{\mathbb R}}\xspace}

\newcommand{\id}{\ensuremath{{\mathbf 1}}\xspace}

\newcommand{\sheaf}[1]{\ensuremath{{\mathcal{#1}}}}

\newcommand{\eqsp}{\quad}
\newcommand{\setsp}{\;\;}

\theoremstyle{plain}
\newtheorem{theorem}{Theorem}[section]
\newtheorem{prop}[theorem]{Proposition}
\newtheorem{lemma}[theorem]{Lemma}
\newtheorem{cor}[theorem]{Corollary}

\theoremstyle{definition}
\newtheorem{definition}{Definition}
\newtheorem{remark}[theorem]{Remark}
\newtheorem{example}[theorem]{Example}

\numberwithin{equation}{section}

\newcommand{\shift}[1]{{#1}^\prime}

\DeclareMathOperator{\spec}{spec}
\DeclareMathOperator{\Span}{Span}
\DeclareMathOperator{\im}{Im}
\DeclareMathOperator{\mes}{mes}

\newcommand{\dmes}[1]{\ensuremath{\,{\mathrm d}#1}}
\newcommand{\dmesns}[1]{\ensuremath{{\mathrm d}#1}} 

\author{Yu. I. Lyubich}
\title{Axiomatic theory of divergent series and cohomological equations}

\date{}

\begin{document}


\bibliographystyle{amsplain}

\maketitle


\begin{abstract}
  A general theory of summation of divergent series based on the
  Hardy-Kolmogorov axioms is developed. A class of functional series
  is investigated by means of ergodic theory. The results are
  formulated in terms of solvability of some cohomological
  equations, all solutions to which are nonmeasurable. In particular, 
  this realizes a construction of a nonmeasurable function as first 
  conjectured by Kolmogorov.
\end{abstract}

\section{Introduction and general theorems}
\label{sec:intro}

\setcounter{equation}{-1}

A natural axiomatic framework for the summation of divergent series already 
appeared in Hardy's early papers, see \cite{hardy74}, and also in
Kolmogorov's short note \cite{kolmogorov25}. Hardy reproduced the
axioms in the book \cite{hardy49} (Section 1.3) and stated that most of the known
summation methods meet them. For instance, this relates to the classical 
Cesaro method $(C,k)$ of any order $k$.

In \cite{lyubich92} the Hardy-Kolmogorov axioms were translated into
the language of functional analysis, and a brief sketch of their main
consequences was presented without proofs. Now we give a developed 
exposition with applications to some functional series generated by dynamical
systems. The ``sums'' of such series satisfy some functional (cohomological) 
equations and, for this reason, they happen to be nonmeasurable.  
Under summability assumption, this phenomenon for a lacunar trigonometric 
series was discovered by Kolmogorov \cite{kolmogorov25} and justified by 
Zygmund \cite{zygmund59}. We prove that functional series in a wide class, 
including Kolmogorov's, are indeed summable. Namely, 
by using the Birkhoff-Khinchin ergodic theorem, we verify that our general
criterion of summability (Theorem \ref{thm:5}) is applicable. Thus,
summation of divergent series is a source of nonmeasureable
functions, as conjectured by Kolmogorov \cite{kolmogorov25}.
 
The following general definition of a summation method is that of \cite{lyubich92}. 
  \begin{definition}
  \label{def:1.0}
  Let $s$ be the linear space of all  sequences $x=(\xi_n\in\C:n\geq 0)$, 
and let $T$ be the shift operator on $s$, i.e. $T(\xi_n) = (\xi_{n+1})$ , 
and, finally, let $L\subset s$ be a $T$-invariant subspace. A linear
  functional $\sigma:L\rightarrow\C$ is called a \defin{summation
    method} (a \defin{summation}, for short) on $L$ if the relation
  \begin{equation}
    \label{eq:1}
    \sigma(x) = \xi_0 + \sigma(Tx),\eqsp x\in L,
  \end{equation}
  is valid. ( In fact, we identify the sequence
  $(\xi_n)$ with the corresponding series $\xi_0 + \xi_1 + \ldots$).
\end{definition}
We do not assume that $s$ is provided with a topology.
  
If there exists a summation $\sigma$ on a $T$-invariant subspace $L$ then we say that
$L$ \defin{admits summation}, and we call $L$ the \defin{domain} of
$\sigma$. In this case for any $T$-invariant subspace $M\subset L$ 
the restriction $\sigma|M$ is a summation on $M$. Moreover, if 
$L= M\oplus N$ where $N$ is also a $T$-invariant subspace then $L$
admits summation if and only if there are some 
summations on $M$ and $N$.

A series $x$ is called \defin{summable} if it belongs to a subspace $L$ 
admitting summation. If the 
summation is $\sigma$, we say that $x$ is \defin{$\sigma$-summable}.

It is very instructive to rewrite (\ref{eq:1}) in the ``cohomological'' form
\begin{equation}
  \label{eq:2}
  \sigma(x) - \sigma(Tx) = \xi_0(x),\eqsp x\in L.
\end{equation}
By iteration of $T$ we obtain
\begin{equation}
  \label{eq:5}
  \sigma(x)-\sigma(T^n x) =s_n(x)\equiv \sum_{k=0}^{n-1}\xi_k ,\eqsp n\geq0.
\end{equation}
(To include the case $n=0$ we set  $s_0(x) = 0$.) By linearity of $\sigma$ 
the formula \eqref{eq:5} can be extended to 

\begin{equation} 
\label{eq:sig}
\phi(1)\sigma(x) - \sigma(\phi(T)x) =  \sum_{n=1}^{m}a_n s_n(x) 
\end{equation}
where  
\begin{equation*}
\phi(\lambda) =  \sum_{n=0}^{m}a_n \lambda^n. 
\end{equation*}
 
A series $x$ is called \defin{finite of length $l(x)=l$} if either $l=0$ (i.e., 
 $x=0$) or $\xi_n =0$ for $n\geq l > 0$, but $\xi_{l-1}\neq 0$.  
The set $F_m$ of finite series $x$ of length $l(x)\leq m$ is a $T$-invariant linear subspace. 
  From \eqref{eq:5} it follows that the functional
  \begin{equation}
    \label{eq:9}
    \sigma_F(x) = \sum_{k=0}^{l(x)-1}\xi_k
  \end{equation}
  \defin{is a unique summation on the space $F$ of all finite series}.

On the space $c^0$ of convergent series we have the \defin{standard summation} 
\begin{equation}
  \label{eq:10}
  \sigma^0(x) = \sum_{k=0}^\infty \xi_k =
  \lim_{n\rightarrow\infty}s_n(x).
\end{equation}
However, \defin{there are infinitely many other summations on $c^0$}. 
We show this after a short considerartion of the general uniqueness problem.

Now we rewrite \eqref{eq:2} as $\sigma(\delta x) = \xi_0(x)$, where  
$\delta = \id - T$ and $\id$ is the identity operator. This $\delta $
is the classical difference operator: 
$\delta(\xi_n) = (\xi_n - \xi_{n+1})$. For every $T$-
invariant subspace $L$ we introduce its \defin{derivative subspace}  
$\shift{L} = \im\delta_L$, where $\delta_L:L\rightarrow L$ is the restriction of $\delta$ to $L$. 
Obviously, $\shift{L}$ is also $T$-invariant.  
 
\begin{lemma}
  \label{prop:2}
  If $L$ admits a summation $\sigma_0$
  then the set of summations on $L$ consists of all linear 
  extensions  of $\sigma_0|\shift{L}$ to $L$. 
\end{lemma}

\begin{proof}
  The equation \eqref{eq:2} is equivalent to
  $\sigma(\delta x) = \sigma_0(\delta x)$, i.e. to $\sigma|\shift{L}=\sigma_0|\shift{L}$.  
\end{proof}

As a consequence, we obtain

\begin{theorem}    
  \label{thm:5.11} 
  A summation on $L$ is unique if and only if $\shift{L} = L$, i.e. the operator
$\delta_L$ is surjective, in other words, the equation $\delta x = y$ has a solution
$x\in L$ for every $y\in L$.   
\end{theorem}

\begin{remark}
\label{rem:0}
\defin{In the whole space} $s$ \defin{the operator} $\delta$ \defin{is surjective, i.e.} 
$\shift{s}= s$. Indeed, for $x = (\xi_n)$ and $y =(\eta_n)$ 
the equation $\delta x = y$ is actually 
$\xi_n - \xi_{n+1} =\eta_n $. Its general solution is $\xi_n =\xi_0 - s_n(y)$
with an arbitrary $\xi_0$, like indefinite integral. 
\end{remark}

\begin{remark}
\label{rem:00}
Using Lemma \ref{prop:2} one can explicitly describe all summations $\sigma$ on $L$.
Namely, for a fixed direct decomposition $L = \shift{L}\oplus R$ we have
$\sigma = \sigma_0\oplus\chi$, where $\chi$ is an arbitraty linear functional on $R$. 
The independent parameters of this description are the values of $\chi$ 
on a basis $B$ of the subspace $R$. We get a one-to-one correspondence 
between summations $\sigma$ on $L$ and complex-valued functions on $B$. As a result,
\defin{if a summation on $L$ is not unique then the set of all summations 
on $L$ is infinite}.
\end{remark}
 
Returning to the space $c^0$ of convergent series we consider 
the closely related space  
$c_0 = \{(\xi_n): \lim_{n\rightarrow\infty} \xi_n =0 \}$ and prove 
\begin{lemma}  
\label{lem:30}
$\shift{c_0} = c^0$, moreover, for every $y =(\eta_n)\in c^0$ 
its unique $\delta$-preimage in $c_0$ is 
\begin{equation*} 
\hat{y} = \left(\sum_{k=n}^\infty \eta_k\right).
\end{equation*}
\end{lemma}

\begin{proof}
If $x = (\xi_n)\in c_0$ then $y=\delta x$ belongs to $c^0$ since 
$s_n(y) =\xi_0 - \xi_n\rightarrow\infty$ as $n\rightarrow\infty$. 
Conversely, if $y\in c^0$ and $x=\hat{y}$ 
then $x\in c_0$ and $\delta x = y$. 
This $x$ is unique since if $\delta x =0$ then all $\xi_n =\xi_0$, hence $x=0$ 
by passing to the limit as $n\rightarrow\infty$.    
\end{proof}

\begin{cor}
\label{cor:0}
  The derivative subspace $\shift{(c^0)}$ is $c^{00} =\{y\in c^0: \hat{y}\in c^0\}$. 
\end{cor}

\begin{proof}
Let $y\in\shift{(c^0)}$, i.e $y=\delta x$, $x\in c^0$. Since $\shift{(c^0)}\subset c^0\subset c_0$,   
we have $y\in c^0$ and $x\in c_0$. By Lemma \ref{lem:30} $\hat{y}=x$, thus $y\in c^{00}$. 
Conversely, let $y\in c^{00}$, i.e.  $y\in c^{0}$ and $\hat{y}\in c^{0}$. Since 
$y=\delta \hat{y}$, we have  $y\in\shift{(c^0)}$.
\end{proof} 

The existence of nonstandard summations on $c^0$ follows from  Theorem \ref{thm:5.11} 
and Corollary \ref{cor:0}. Indeed, the set $c^0\setminus c^{00}$ is not empty. For instance, 
it contains any series
\begin{equation*}
  \zeta_{\alpha} = \left((n+1)^{-\alpha}\right), \eqsp 1<\alpha\leq 2.
\end{equation*}
  
Now we proceed to the general existence problem. 

\begin{lemma}
\label{lem:1.7}
The series $\pi_0: 1+1+\ldots$ is not summable.
\end {lemma}

\begin{proof} 
Let $\pi_0$ be $\sigma$-summable. Then $\sigma(\pi_0)- \sigma(T\pi_0) = 0$
since $T\pi_0 = \pi_0$. On the other hand, $\xi(\pi_0)=1$. 
\end{proof} 

The series $\pi_0$ generates the $1$-dimensional subspace $\Pi_0 =\ker\delta$ 
of {\em constant series}, $\xi_0 + \xi_0 + \ldots$. This is a 
subspace of the space $\Pi_\infty$ of the series
$\left(\pi(n)\right)$, where $\pi$ runs over all 
complex-valued polynomials. Obviously, 
$\Pi_\infty$ is $T$-invariant as well as 
every its subspace $\Pi_m = \{\pi: \deg\pi\leq m\}$, $m\geq 0$. 
Moreover,  $\shift{\Pi_m} = \Pi_{m-1}, m\geq 1$.  

\begin{theorem}
  \label{thm:5}
  Let $L$ be a $T$-invariant subspace. The following
  statements are equivalent.
  \begin{enumerate}
    \item $L$ admits summation.
    \item $\pi_0$ does not belong to $L$, i.e. $L\cap\Pi_0 =0$
    \item There is no nonzero polynomial series in $L$, i.e. 
$L\cap\Pi_\infty = 0$. 
\end{enumerate}
\end{theorem}
\begin{proof}
  $(1)\Rightarrow(2)$ since $\pi_0$ is not summable. Conversely,
  $(2)\Rightarrow(1)$. Indeed, we have $L\cap\ker\delta = 0$, hence
  $\delta_L$ is injective, so left invertible. 
  Let $i:L\rightarrow L$ be a left inverse to $\delta_L$. 
  (This is not unique if $\shift{L}\neq L$.)
  Then the linear functional $\sigma_0(x) = \xi_0(ix)$ is a summation 
  on $L$ since $\sigma_0(\delta x) = \xi_0(x)$.

Obviously, $(3)\Rightarrow(2)$. 
Conversely, $(2)\Rightarrow(3)$. Indeed, let $\pi\in L\cap\Pi_{\infty}$, 
and let $\pi\neq 0, \deg\pi =m$. Then $\delta^m \pi = \gamma\pi_0$ where 
$\gamma$ = const $\neq 0$. Hence $\pi_0\in L$ in contrary to (2).
\end{proof}

\begin{cor}
  \label{thm:6.11}
  A $T$-invariant subspace $L$ admits summation if and only
  if the operator $\delta_L$ is injective.  
\end{cor}

Combining this result with Theorem \ref{thm:5.11} we obtain

\begin{cor} 
  \label{thm:6.12} 
  A $T$-invariant subspace $L$ admits a unique summation if and only 
  if the operator $\delta_L$ is bijective.
\end{cor} 

\begin{cor}
   \label{cor:6}
  If a finite-dimensional $T$-invariant subspace admits summation 
then the summation is unique.
\end{cor}

\begin{proof}
In this case all injective operators are surjective. 
\end{proof}

\begin{remark} 
\label{rem:01}
Informally speaking, the subspaces with a unique summation are 
just those where the "integration" becomes definite. Indeed, 
$\delta_L$ is bijective if and only if it is invertible. 
\end{remark} 

\begin{remark}
\label{rem:02}
For any $T$-invariant subspace $L$ let us consider the sequence 
of derivative subspaces 
$L\supset\shift{L}\supset\shift{\shift{L}}\supset\ldots$. Let $N$ 
be their intersection. One can prove that \defin{if $L$ admits summation
but this is not unique then this remains not unique on all derivative subspaces, 
while becomes unique on $N$.} 
\end{remark}

The summability problem can be "localized" as follows. 

\begin{theorem}  
\label{thm:6.13} 
Let $L$ be a $T$-invariant subspace. Then $L$ admits summation if and only if
every $x\in L$ is summable.
\end{theorem}

\begin{proof} 
The "only if" part is trivial. The "if" part follows from Theorem \ref{thm:5} since
$\pi_0$ is not summable, thus $\pi_0\not\in L$.
\end{proof}

Now for every $x\in s$ we consider the smallest $T$-invariant 
subspace $L_x$ containing $x$. This is  
\begin{equation*}
  L_x = \Span(T^n x)= \{\phi(T)x\}  
\end{equation*} 
where $\phi$ runs over all polynomials of one variable. Obviously, 
$x$ is summable if and only if $L_x$ admits summation. Combining this fact 
with Theorem \ref{thm:6.13} we obtain a "local" version of Theorem \ref{thm:5}. 

\begin{theorem} 
\label{thm:2loc}
A $T$-invariant subspace $L$ admits summation if and only if for every
$x\in L$ the subspace $L_x$ does not contain $\pi_0$ or, equivalently,
$L_x\cap\Pi_\infty = 0$.
\end{theorem}  

We conclude this section with a few examples.

\begin{example}
  \label{ex:1}
     The subspace $c_0$ admits summation since $\pi_0\not\in c_0$.
\end{example}

\begin{example}
  \label{ex:2}
 The subspace $m = \{x:\sup_n |\xi_n| < \infty \}$ does not admit summation 
 since $\pi_0\in m$. However, given a Banach limit (an invariant mean) on $m$, 
 the $T$-invariant 
 subspace $m_0 = \{x \in m: B$-$\lim_{n\rightarrow\infty} \xi_n =0 \}$ 
 admits summation since $\pi_0\not\in m_0$. Note that $m_0\supset c_0$ since 
 the Banach limit coincides with the standard limit for all convergent sequences. 
\end{example}

\begin{example}
  \label{ex:3}
 The formula $\sigma(x) = B$-$\lim_{n\rightarrow\infty}s_n(x)$
 determines a summation on the subspace  $\hat{m} = \{x:\sup_n |s_n(x)| < \infty \}$. 
 We call it the \defin{Banach summation}. 
 Note that $\hat{m}\subset m_0$ since $\xi_n = s_{n+1}- s_n$ and the Banach limit is $T$-invariant.
\end{example}

\begin{example}
  \label{ex:4}
    A summation method going back to Euler uses an analytic continuation 
of the \defin{generating function}
  \begin{equation}
    \label{eq:2.12}
     g(t;x) = \sum_{n=0}^\infty \xi_n t^n, \eqsp x=(\xi_n),
  \end{equation}
  where $t$ is a complex variable. This function is defined and analytic
  in the disk $|t|<r_x$ if 
  \begin{equation}
    \label{eq:2.13}
    r_x \equiv \left( \limsup_{n\rightarrow\infty}
      |\xi_n|^{1/n}\right)^{-1} > 0.
  \end{equation}
Assume that all series from a $T$-invariant subspace $L$
  satisfy (\ref{eq:2.13}), so we have a linear space $\sheaf{A}_L$ of
  analytic germs $g(t;x)$ at $t=0$. Let $G\subset \C$ be an open connected 
  set containing $t=0,1$, and let each $g\in\sheaf{A}_L$ be the Taylor germ 
  of a function $\tilde{g}(t;x)$ that is analytic on $G_x =G\setminus\Gamma_x$ 
  where $\Gamma_x$ is a finite set, $1\not\in\Gamma_x$. Then 
$\epsilon(x)= \tilde{g}(1;x)$ is a summation on $L$ (\defin{Euler's summation}).
 Indeed, from (\ref{eq:2.12}) it follows that
  \begin{equation*}
    g(t;x) - tg(t;Tx) = \xi_0 (x),\eqsp |t|<\min(r_x, r_{Tx}),
  \end{equation*}
  By uniqueness of the analytic continuation $\tilde{g}(t;x)$ 
is a linear functional of $x$ and
  \begin{equation*}
    \tilde{g}(t;x) - t\tilde{g}(t;Tx) = \xi_0(x),\eqsp t\in G_x\cap G_{Tx}.
  \end{equation*}
  In particular, $ \epsilon(x)-\epsilon(Tx)=\tilde{g}(1;x) - \tilde{g}(1;Tx) = \xi_0(x)$.
\end{example}

The best known example of Euler's sum is 
$\epsilon (((-1)^n)) = 1/2$. More generally,
\begin{equation*}
\epsilon\left( (\lambda^n)\right) = 
(1-\lambda)^{-1}, \eqsp \lambda \in \C\setminus\{1\}.
\end{equation*}

\section{Quasiexponential series}
\label{sec:2}

The geometric progression (or \defin{exponential series}) $(\lambda^n)$ 
is an eigenvector 
of $T$ for the eigenvalue $\lambda \in \C$ . (For $\lambda=0$ we set $0^0=1$.)
The corresponding eigenspace $E_\lambda = \ker(T-\lambda \id)$ is $1$-dimensional.
This is the first member of the increasing sequence of the \defin{root subspaces} 
$E_{\lambda,m} = \ker(T-\lambda \id)^m$,
$m=1,2,3,\ldots$; their union is denoted by $E_{\lambda,\infty}$. We have $E_{0,m} = F_m$, the 
space of all finite series of length $\leq m$, so $E_{0,\infty} = F$.
If $\lambda\neq0$ then $E_{\lambda,m}$ consists of
all series $(\pi(n)\lambda^n)_0^\infty$, $\pi\in\Pi_{m-1}$. 
In particular, $E_{1,m} = \Pi_{m-1}$, $E_{1,\infty} = \Pi_{\infty}$. 
Note that $\dim E_{\lambda,m} = m$ in any case.

Now let $\phi(\lambda)$ be a nonconstant polynomial, i.e.
\begin{equation*}
  \phi(\lambda) = \lambda^{m_0} \prod_{k=1}^\nu
  (\lambda-\lambda_k)^{m_k},  
\end{equation*}
where $m_0\geq0$, and if $\nu>0$ then $\lambda_k$ are nonzero 
pairwise distinct roots with multiplicities $m_k\geq1$. Then
\begin{equation*}
  \ker\phi(T) = E_{0,m_0} \oplus E_{\lambda_1,m_1}
  \oplus\cdots\oplus E_{\lambda_\nu,m_\nu}
\end{equation*}
according to the Jordan form of $T|\ker\phi(T)$.
In other words, $\phi(T)x = 0$ \defin{if and only if} 
\begin{equation}
    \label{eq:2.6}
    \xi_n = \sum_{k : \lambda_k\neq 0} \pi_k(n){\lambda_k^n} + \zeta_n,
\end{equation}
  \defin{where} $\pi_k\in\Pi_{m_k-1}$ \defin{and} $(\zeta_n)\in F$,
  $\zeta_n=0$ \defin{for} $n\geq m_0$. The decomposition \eqref{eq:2.6} is unique.
By the way, $\phi(T)x = 0$ is nothing but a homogeneous linear difference 
equation with constant coeeficients, and \eqref{eq:2.6} is its general solution. 
In the case $m_0 =0$ this formula turns into the classical one 
concerning the two-sided sequences. 

Any series $x$ with the members of form \eqref{eq:2.6} is called \defin{quasiexponential}. 
The complex linear space of all quasiexponential series will be denoted by $Q$. This is a
$T$-invariant subspace of $s$. 
Letting $Z_{m_0}=\{0\}$ for 
$m_0>0$ and $Z_{m_0}=\emptyset$ for $m_0=0$, we call 
the set $\{\lambda_k : \pi_k\neq 0\} \cup Z_{m_0}$ the 
\defin{spectrum} of $x\in Q$ and denote it by $\spec(x)$. Obviously, 
$\spec(x)\neq\emptyset$ if $x\neq 0$ but $\spec(0)=\emptyset$. 
Let $x\neq 0$. Then $x$ is finite or polynomial
if and only if $\spec(x) = \{0\}$ or $\{1\}$, respectively. 
In general, the $\spec(x)$ coincides with the set of roots of a 
minimal polynomial $\phi_x(\lambda)$ such that $\phi_x(T)x=0$.
(As usual, the minimality means that $\deg\phi_x$ is minimal. 
This polynomial is unique up to a constant factor.) 

The following theorems show the importance of the quasiexponential 
series for the general summation theory. 
 
\begin{theorem}
  \label{thm:2.1}
  A series $x$ is not summable if and only if $x\in Q$ and $1\in\spec(x)$.  
\end{theorem}
\begin{proof}
By Theorem \ref{thm:2loc} $x$ is not summable if and only if $\pi_0\in L_x$, i.e.
$L_x\cap\ker(\id - T)\neq 0$ or, equivalently, there is 
a polynomial $\psi$ such that 
\begin{equation}
\label{eq:A} 
\psi(T)x\neq 0, \eqsp (\id -T)\psi(T)x = 0. 
\end{equation}
From \eqref{eq:A} it follows that $x\in Q$ and $(1-\lambda)\psi(\lambda)$ 
is divided by the minimal polynomial $\phi_x(\lambda)$. Hence, $\phi_x(1) = 0$, otherwise, 
$\psi$ is divided by $\phi_x$, so $\psi(T)x = 0$. Conversely, if $x\in Q$ and $\phi_x(1)=0$ 
then  $\phi_x(\lambda) = (1-\lambda)\psi(\lambda)$, so \eqref{eq:A} is valid 
since $\deg\psi<\deg\phi_x$.      
\end{proof}

In view of Theorem \ref{thm:2.1} let us introduce 
\begin{equation*}
Q_1 = \{x\in Q: 1\not\in\spec(x)\}, 
\end{equation*} 
so $x\in Q_1$ if and only if $\phi(T)x = 0$ for a polynomial $\phi$ 
such that $\phi(1)\neq 0$. The subspace $Q_1$ is $T$-invariant, and  
$Q = Q_1 \oplus\Pi_\infty$ according to \eqref{eq:2.6}. 

Now note that $x\in Q$ if and only if 
the set $\{T^n x\}_0^\infty$ is linearly dependent, i.e the subspace $L_x$ 
is finite-dimensional. Its basis is $\{T^n x\}_0^{\nu-1}$, where $\nu$ is the degree
of the related minimal polynomial, thus, $\dim L_x=\nu$. If $x\not\in Q$
then $L_x$ is iinfinite-dimensional since $\{T^n x\}_0^\infty$is its basis.

\begin{theorem}
  \label{thm:2.2}
Let a series $x$ be summable. Then the summation on $L_x$ is unique 
if and only if $x\in Q_1$ .
\end{theorem}
\begin{proof}
By Theorem \ref{thm:2.1} either $x\in Q_1$,   
or $x\not\in Q$. In the first case the summation is unique by Corollary \ref{cor:6}.
In the second case the values $\sigma(T^n x), n\geq 1$, are determined by \eqref{eq:5}, 
while $\sigma(x)$ remains arbitrary.
\end{proof}

\begin{cor}
\label{cor:2.10}
Every subspace $L\subset Q_1$ admits a unique summation.
\end{cor}

\begin{proof}
Let $x\in L$ and let $\sigma$ be a summation on $L$. Then $\sigma(x) = (\sigma|L_x)(x) $.
By Theorem \ref{thm:2.2} the summation $\sigma|L_x$ does not depend on choice of $\sigma$.
\end{proof}

We denote by $\epsilon_1$ the unique summation on $Q_1$. 
By Corollary \ref{cor:2.10} the unique summation on any subspace 
$L\subset Q_1$ is $\epsilon_1|L$. 

\begin{example}
\label{ex:2.0}
Let $Q_0 = Q\cap{c^0}$, i.e. $Q_0$ is the subspace of convergent quasiexponential series.
Using \eqref{eq:2.6} one can prove that $x\in Q_0$ if and only if 
$x\in Q$ and $\spec(x)\subset \{\lambda\in \C: |\lambda|<1\}$. 
In particular, $1\not\in \spec(x)$ for $x\in Q_0$, so $Q_0 = Q_1\cap{c^0}$. Therefore, 
on $Q_0$ the summation $\epsilon_1$ coincides with the standard summation $\sigma^0$.
\end{example}

Now we prove that \defin{the  summation $\epsilon_1$ coincides with the 
restriction of Euler's summation $\epsilon$ to the subspace $Q_1$}.   
\begin{lemma}
  \label{lem:2.5}
  For $x\in Q$ the generating function $g(t,x)$ is well-defined, 
and $\tilde{g}(t;x)$ is a rational function of $t$. The set
  of its poles is
  \begin{equation}
    \label{eq:2.15}
    \left\{ t = \lambda^{-1} : \lambda\in\spec(x),\setsp \lambda\neq 0\right\}.
  \end{equation}
\end{lemma}
\begin{proof}
  It suffices to consider the case of a single-point spectrum. If
  $\spec(x) = \{0\}$ then $x$ is finite, so $\tilde{g}(t;x)$ is a
  polynomial in $t$. On the other hand, 
the set \eqref{eq:2.15} is empty in this case. Now let $\spec(x)=\{\lambda\}$, 
  $\lambda\neq 0$. Then $x = \left(\pi(n)\lambda^n\right)$ where
  $\pi\in\Pi_\infty$. Accordingly,
  \begin{equation*}
    g(t;x) = \sum_{n=0}^\infty \pi(n)\lambda^n t^n, \eqsp |t|<|\lambda|^{-1}.
  \end{equation*}
  If $\deg\pi=\nu-1$ then $\pi$ can be represented as
  \begin{equation*}
    \pi(n) = \sum_{k=0}^{\nu-1} c_k\binom{k+n}{k}, \eqsp c_{\nu-1}\neq 0.
  \end{equation*}
  This yields
  \begin{equation}
    \label{eq:2.16}
    \tilde g(t;x) = \sum_{k=0}^{\nu-1}\frac{c_k}{(1-\lambda t)^{k+1}}.
  \end{equation}
\end{proof}

By the way, any rational function $g(t)$ which is regular at $t=0$ is 
the generating function of a series $x\in Q$. According to
  (\ref{eq:2.16}), this $x$ can be obtained from the decomposition of
  $g(t)$ into partial fractions.
 
\begin{cor}
\label{thm:2.6}
  For $x\in Q_1$ the function $\tilde g(t;x)$ is rational and regular at $t=1$.  
\end{cor}

This means that $Q_1$ is a subspace of the domain of Euler's summation. 
Therefore, $\epsilon_1 = \epsilon |Q_1$ by uniqueness of summation on $Q_1$.
   
For $x\in Q_1$ an explicit expression of $\epsilon_1 |L_x$ 
follows from the formula \eqref{eq:sig}. Namely, \defin{if} 
\begin{equation*}
\phi_x(\lambda) = \sum_{n=0}^{\nu}a_n\lambda^n
\end{equation*}
\defin{is a corresponding minimal polynomial then} 
 
\begin{equation*}
\epsilon_1(z) = \frac{1}{\phi_x (1)}\sum_{n=1}^{\nu}a_n s_n(z), \eqsp z\in L_x . 
\end{equation*}
Indeed, $\phi_x(T)z = 0$ for all $z\in L_x$. It remains to 
substitute  $x=z$ and $\phi =\phi_x$ into \eqref{eq:sig}. 
Actually, we see that the minimal polynomial $\phi_x$ can be changed 
to any polynomial $\phi$ such that $\phi(T)x =0$ and $\phi(1)\neq 0$.

\begin{example}
\label{ex:2.1}
Let $x$ be $(l+1)$-periodic, i.e. $T^{l+1}x = x$, and let $1\not\in\spec(x)$, 
or, equivalently, $s_{l+1}(x)= 0$. Then
\begin{equation*}
\epsilon_1(x) =\frac{1}{l+1} \sum_{n=1}^{l}s_n(x),
\end{equation*}
so $\epsilon_1(x)$ coincides with the Cesaro sum of order 1. 
This is not an occasional fact. The point is that     
a quasipolynomial series $x$ is $(C,1)$-summable if and only if  
$\spec(x)\subset \{\lambda\in \C: |\lambda|\leq 1\}$ and the  
roots of the minimal polynomial luying on the unit circle 
are simple. If, in addition, 
$1\not\in\spec(x)$ then the Cesaro sum of $x$ 
coincides with $\epsilon_1(x)$ by uniqueness. 
\end{example}

\section{Extension theory}
\label{sec:3}

In spirit of the classical definition (see e.g. \cite{hardy49}, Section 4.3) 
we say that a summation method
$\tau$ is \defin{stronger} than $\sigma$ and write $\tau\succ\sigma$ if 
$\tau$ is an extension of $\sigma$.  
For instance, $\sigma^0\succ\sigma_F$, see \eqref{eq:9} and\eqref{eq:10}.
In turn,   
a method $\sigma$ is called \defin{regular} if $\sigma\succ\sigma^0$.
The Banach summation on $\hat{m}$ (Example \ref{ex:3}) is regular by definition of the 
Banach limit. The Cesaro methods of all orders are regular,
while Euler's method is not. Indeed, though $r_x\geq1$ for  $x \in c^0$ 
 but for some $x$'s the function $g(t;x)$ cannot be analytically continued 
to $t=1$. 

A summation $\mu$ is called \defin{maximal} if there is no 
summations $\tau\succ\mu, \tau\neq\mu$. 

\begin{theorem}
\label{thm:3.3}
For every summation $\sigma$ there exists a maximal $\mu\succ\sigma$.
\end{theorem}
\begin{proof}
We use the Zorn lemma. The relation ``$\succ$'' is a partial order 
on the set of all summations, a fortiori, on 
the subset $\{\tau:\tau\succ\sigma\}$.  
This order is inductive: there is a majorant 
$\tau$ for any linearly ordered subset $\{\tau^{\alpha} \succ\sigma\}$. 
Indeed, let $L^{\alpha}$ be the domain of $\tau^{\alpha}$, and let  
$ L = \bigcup_\alpha L^{\alpha}$. Then $\tau$ is
well-defined on $L$ as $\tau x = \tau^{\alpha_x} x$ where 
$\alpha_x$ is any index such that $x\in L^{\alpha_x}$.
\end{proof}

\begin{cor}
\label{cor:3.0}
There exists a regular maximal summation.
\end{cor}

Actually, any extension can be realized as a sequence (transfinite, in general) 
of minimal steps. Every such step extends the domain $L$ of a summation 
$\sigma$ to $L[x] = L+L_x$ with some $x\not\in L$. 
To analyze this  situation we consider the set $I_{x,L}$ of polynomials 
$\phi(\lambda)$ such that $\phi(T)x\in L$. Since $L$ is $T$-invariant,  
the $I_{x,L}$ is an ideal of the ring of all polynomials of $\lambda$. 
Though $0\in I_{x,L}$, the nonzero constants 
do not belong to $I_{x,L}$ as long as $x\not\in L$. 
Obviously, $I_{x,K}\subset I_{x,L}$ if $K\subset L$, in particular, 
$I_{x,0}\subset I_{x,L}$. Implicitly, we already dealt with the ideal $I_{x,0}$ 
in Section \ref{sec:2}.

\begin{lemma}
\label{lemma:3.2} 
If $I_{x,L} = 0$ then any summation $\sigma$ on $L$ extends to $L[x]$.
\end{lemma}

\begin{proof}
In this case $L\cap L_x = 0$ and $x\not\in Q$. Thus, $L[x] = L\oplus L_x$ 
and $L_x$ admits summation by Theorem \ref{thm:2.1}. 
\end{proof}

\begin{remark}
In Lemma \ref{lemma:3.2} the set of extensions is infinite by Theorem \ref{thm:2.2}.
\end{remark}

Now we assume $I_{x,L} \neq 0$  and introduce
\begin{equation*} 
\theta _{x,L}(\lambda) = \sum_{n=0}^{\nu}c_n\lambda ^n , 
\end{equation*}
a minimal polynomial in $I_{x,L}$. This is a greatest common divisor of all $\phi\in I_{x,L}$. 
Below we use the reduced notation $\theta(\lambda) \equiv\theta_{x,L}(\lambda)$.
It is convenient to normalize this polynomial so that $c_{\nu} = 1$, i.e.
\begin{equation}
\label{eq:3.11}
\theta(\lambda) = \lambda^{\nu} +  \sum_{n=0}^{\nu - 1}c_n\lambda ^n ,
\end{equation}
The trivial case $\theta(\lambda) \equiv 1$ (i.e. $\nu = 0$) is formally included
in this setting.
  
\begin{lemma} 
\label{lemma:3.3}
A summation $\sigma$ on $L$ extends to $L[x]$ if and only if 
either $\theta(1)\neq 0$, and then $\sigma$ is arbitrary, 
or $\theta(1)=0$, and then $\sigma$ is such that
\begin{equation} 
\label{eq:3.1}
\sigma (\theta(T)x )+ \sum_{n=1}^{\nu}c_n s_n(x) = 0.
\end{equation}
The extension is unique if and only if $\theta(1)\neq 0$.
\end{lemma}
 
\begin{proof}
Let $\tau$ be an extension of $\sigma$ to $L[x]$. Then
\begin{equation} 
\label{eq:3.2}
\theta(1)\tau(x) = \tau(\theta(T)x )+ \sum_{n=1}^{\nu}c_n s_n(x)
\end{equation} 
according to \eqref{eq:sig}. However, $\tau(\theta(T)x)=\sigma (\theta(T)x)$ since 
$\theta(T)x\in L$ and $\tau|L =\sigma$. Thus, if $\theta(1)=0$ then \eqref{eq:3.2}
turns into \eqref{eq:3.1}.
 
In the converse direction we start with a value $\tau(x)$ such that
\begin{equation} 
\label{eq:3.3}
\theta(1)\tau(x) = \sigma(\theta(T)x )+ \sum_{n=1}^{\nu}c_n s_n(x).
\end{equation}
This value does exist under our conditions (and uniquely determined 
if $\theta(1)\neq 0$, otherwise, it is arbitrary). Setting 
\begin{equation}
\label{eq:3.4}  
\tau(T^n x) = \tau(x) - s_n(x), \eqsp 1\leq n\leq\nu-1 ,
\end{equation}
we determine a linear extension $\tau$ of $\sigma$ to 
\begin{equation*}
L[x]=L\oplus R, \eqsp R=\Span\{T^n x\}_0^{\nu-1}.
\end{equation*}
To prove that $\tau$ is a summation it remains to verify the equality 
\begin{equation}
\label{eq:3.5}
\tau(T^{\nu} x) = \tau(x) - s_{\nu}(x).
\end{equation}

Note that, as a rule, $T^{\nu}x\not\in R$, so the space $R$ is not $T$-invariant. 
Indeed, by \eqref{eq:3.11} we have  
\begin{equation}
\label{eq:3.6}
T^{\nu}x = \theta(T)x\oplus(T^{\nu} -\theta(T))x  ,
\end{equation}
so $T^{\nu}x\not\in R$ as long as $\theta(T)x\neq  0$. According to \eqref{eq:3.6}, 
\begin{equation*} 
\tau(T^{\nu}x) = \sigma(\theta(T)x) + \tau((T^{\nu} -\theta(T))x) = 
\sigma(\theta(T)x) - \sum_{n=0}^{\nu-1}c_n\tau(T^n x).  
\end{equation*} 
By substitution from \eqref{eq:3.4} we obtain  
\begin{equation*}
\tau(T^{\nu}x) = \sigma(\theta(T)x) + 
\sum_{n=0}^{\nu-1}c_n s_n(x) -
\tau(x)\sum_{n=0}^{\nu-1}c_n . 
\end{equation*} 
(Recall that $s_0(x)=0$). This yields \eqref{eq:3.5} 
because of \eqref{eq:3.3} and the relation 
\begin{equation*} 
\theta(1) -  \sum_{n=0}^{\nu-1}c_n  = c_\nu =1 .
\end{equation*}    
\end{proof}

\begin{cor} 
\label{cor:3.6} 
Every maximal summation $\mu$ is stronger than $\epsilon_1$.
\end{cor}
 
\begin{proof} 
By the unqueness of the summation $\epsilon_1$ on $Q_1$ 
we only have to show that the domain $M$ of $\mu$ contains $Q_1$. 
Let $x\in Q_1$, so $\phi(T)x=0$ where $\phi$ is a polynomial,
$\phi(1) \neq0$. Thus, $\phi\in I_{x,M}$, so $\theta_{x,M}$ is a 
divisor of $\phi$. Therefore, $\theta_{x,M}(1)\neq 0$. By Lemma
\ref{lemma:3.3} and maximality of $\mu$ we obtain $x\in M$.  
\end{proof}

Lemma \ref{lemma:3.3} shows that the only obstacle to extension of a summation 
$\sigma$ from $L$ to $L[x]$ is the inequality 
\begin{equation*}
\sigma (\theta(T)x )+ \sum_{n=1}^{\nu}c_n s_n(x) \neq 0
\end{equation*}
in the case $\theta(1)=0$. 
However, this  obstacle is removable by a "polynomial regularization"  of $x$. 
 
\begin{lemma}
\label{lemma:3.4}
Let $\theta_{x,L}(1) =0$ and let $m$ be 
the multiplicity of this root of $\theta_{x,L}(\lambda)$.
Then there exists a polynomial series $\pi$ 
of degree $\leq m-1$ such that any summation $\sigma$ 
extends from $L$ to $L[x-\pi]$.
\end{lemma} 

\begin{proof}
We start with the case $\theta_{x,L}(\lambda) = (\lambda -1)^ m$.
For every $y=x-p$, $p\in \Pi_{m_x - 1}$, we have $\theta_{x,L}(T)y =\theta_{x,L}(T)x\in L$,  
so $\theta_{x,L}\in I_{x,L}$. Hence, $\theta_{x,L}$ is divided by $\theta_{y,L}$, so 
$\theta_{y,L}(\lambda) = (\lambda -1)^ {m_y}$ with some $m_y\leq m$. 
We choose the summand $p$ in $y$ to make $m_y$ minimal. 
If $m_y = 0$ then $\theta_{y,L} = 1$, hence $y\in L$. Thus, we have a trivial 
extension $L[x-\pi] = L$ with $\pi = p$. 

Let $m_y\geq 1$.  Then we consider $z=y-q$, $q\in \Pi_{m_y-1}$, 
so that $z = x - \pi$ where $\pi = p+q\in \Pi_{m -1}$. 
As before,  $\theta_{z,L}(\lambda) = (\lambda -1)^ {m_z}$ 
where $m_z\leq m_y$. Finally, $m_z = m_y$ by minimality of the latter. 
Thus, $\theta_{z,L} = \theta_{y,L}$, and, accordingly,
\begin{equation} 
\label{eq:3.a}  
\sigma (\theta_{z,L}(T)x )+ \sum_{n=1}^{m_y}c_n s_n(z) = 
\sigma (\theta_{y,L}(T)y )+ \sum_{n=1}^{m_y}c_n s_n(y) -\sum_{n=1}^{m_y}c_n s_n(q) .
\end{equation} 
The corresponding obstacle to extension of $\sigma$ to $L[z]$ 
disappears if, for instance,   
\begin{equation*}
q(n) = \alpha \binom{n}{m_y - 1}
\end{equation*}
with a suitable $\alpha\in\C$. Indeed, for this $q$ the subtrahend in \eqref{eq:3.a} 
reduces to $\alpha$.   

In general, $\theta_{x,L}(\lambda) = \phi(\lambda)(\lambda -1)^m$, where $\phi(1)\neq 0$.
With $u = (T-1)^mx$ we have $\phi(T)u = \theta_{x,L}(T)x\in L$. Therefore, $\phi$
is divided by $\theta_{u,L}$, thus $\theta_{u,L}(1)\neq 0$. By Lemma \ref{lemma:3.3}
$\sigma$ extends to a summation $\tau$ on $L[u]$. In turn, 
$\theta_{x,L[u]}(\lambda) = (\lambda -1)^ l$ with $l\leq m$. 
Hence, there exists $\pi\in \Pi_{l-1}\subset\Pi_{m-1}$ such that 
$\tau$ extends to $(L[u])[x-\pi]$, and, a fortiori, to $L[x-\pi]$.
\end{proof}
Combining Lemmas \ref{lemma:3.2}, \ref{lemma:3.3} and \ref{lemma:3.4} we obtain 
the following general

\begin{theorem} 
\label{thm:3.6} 
Let a subspace $L$ admit summation. For every $x\in s$ there exists a polynomial
series $\pi$ such that any summation $\sigma$ extends from $L$ to $L[x-\pi]$.   
\end{theorem} 

As an important consequence we obtain

\begin{theorem}
\label{thm:3.66}
A $T$-invariant subspace $M$ is the domain of a maximal summation if and only if
\begin{equation}
\label{eq:3.4p}
M\oplus \Pi_\infty = s,
\end{equation}
i.e. $M$ is a $T$-invariant direct complement of the subspace of the  
polynomial series to the whole space $s$. 
\end{theorem}

Thus, every maximal summation is applicable to all 
series up to a polynomial regularization. In this sense, the maximal 
summations are \defin{universal}. 

\begin{proof}
"If". $M$ admits summation, since $M\cap \Pi_\infty = 0$. Any summation on $M$ 
is maximal since any nontrivial extension of $M$ intersects $\Pi_\infty$.

"Only if". We apply Theorem \ref{thm:3.6} to $L=M$. 
By maximality of $M$ the extension $M[x-\pi]$ is trival, i.e. 
$x-\pi\in M$. Thus, $M + \Pi_\infty = s$.
Moreover, $M\cap\Pi_{\infty} =0$ since $M$ admits summation.   
\end{proof}   

\begin{cor}
\label{cor:3.5}
Every maximal summation is unique on its domain.
\end{cor}

\begin{proof}
The operator $\delta =\id -T$ is surjective on the whole space $s$, 
see Remark \ref{rem:0}.
Since in \eqref{eq:3.4p} both summands are $T$-invariant, the
restriction $\delta_M$ is also surjective.
Thus, Theorem \ref{thm:5.11} is applicable. 
\end{proof}

As a result, we have {\em a 1-1 correspondence between maximal 
summations and $T$-invariant direct complements of $\Pi_\infty$ to $s$}.

\section{Orbital series}
\label{sec:5}

A \defin{functional series on a set} $A\neq\emptyset$ is a
mapping $X:A\rightarrow s$, i.e. for every $\alpha\in A$ we have
a numerical series $X(\alpha) = (\xi_n(\alpha))$. Given a
summation $\sigma$ with a domain $L$, we say that $X$ is
\defin{$\sigma$-summable} if such are all series $X(\alpha)$,
i.e. $\im X\subset L$ and 
\begin{equation}
 \label{eq:5.30}
  \sigma(X(\alpha))-\sigma((TX)(\alpha)) = \xi_0(\alpha),\eqsp \alpha\in A, 
\end{equation}
where $(TX)(\alpha) = T(X(\alpha)) = (\xi_{n+1}(\alpha))$.

A functional series $X$ is called \defin{summable} if there exists a
summation $\sigma$ such that $X(\alpha)$ is $\sigma$-summable
for every $\alpha$. For example, every trigonometric series whose coefficients 
tend to zero is summable. Moreover, there is a common summation for all 
these series, namely, any summation on $c_0$.  

An important class of functional series is 
\begin{equation}
  \label{eq:5.5}
  X(\alpha) = (\xi_0(f^n\alpha)), \eqsp \alpha\in A,
\end{equation}
where $f$ is a mapping $A\rightarrow A$. 
The sequence $(f^n\alpha)$ is the $f$-orbit of the point
$\alpha$, therefore, we call the functional series \eqref{eq:5.5} 
\defin{orbital}. In this case the subspace 
\begin{equation*}
L_X = \Span (\im X)\subset s 
\end{equation*}
is $T$-invariant since
\begin{equation}
\label{eq:5.55}
(TX)(\alpha) = X(f\alpha).
\end{equation} 
Hence, \defin{an orbital series $X$ is summable if and only there exists a 
summation on $L_X$}. Combining \eqref{eq:5.55} and \eqref{eq:5.30} we obtain 

\begin{prop}
 \label{prop:5.1}
If an orbital series $X$ is $\sigma$-summable then the function
$\psi(\alpha) = \sigma(X(\alpha))$ satisfies the cohomological equation (\em{c.e.})
\begin{equation}
  \label{eq:5.8}
  \psi(\alpha) - \psi(f\alpha) = \xi_0(\alpha), \eqsp \alpha\in A.
\end{equation}
\end{prop}
This is a bridge between the summations and the functional equations playing 
a considerable role in the modern theory of dynamical systems 
and group representation theory, see e.g. \cite{anosov73}, \cite{gottschalk55}, 
\cite{katok95}, \cite{kirillov67}. 
In standard terms related to the dynamical system $(A,f)$, any function $\psi : 
A\rightarrow \C$ is a 0-cochain , and its coboundary is the 1-cochain 
\begin{equation*}
\theta(n,\alpha)= \psi(\alpha) - \psi(f^n\alpha), \eqsp n\geq 0,\eqsp \alpha\in A.
\end{equation*}
A 1-cochain $\omega(n,\alpha)$ is a cocycle if 
\begin{equation*}
\omega(n, f^m\alpha)- \omega(n+m,\alpha) + \omega(m,\alpha) = 0 \eqsp (n,m\geq 0). 
\end{equation*}
Every coboundary is a cocycle but, in general, the converse is not true, i.e. not every cocycle 
is "cohomologically trivial". For any 0-cochain $\xi_0$ the 1-cochain 
\begin{equation*}
s(n,\alpha)\equiv s_n(\alpha)=\sum_{k=0}^{n-1}\xi_0(f^k\alpha)   
\end{equation*}
is a cocycle. This cocycle is a coboundary if and only if c.e. \eqref{eq:5.8} is solvable. 

In the context of summations we have a dynamical system $(L,T)$, where $L$ is a $T$-invariant 
subspace of $s$, and deal with the cocycle 
\begin{equation}
\label{eq:coc}
s(n,x)\equiv s_n(x)=\sum_{k=0}^{n-1}\xi_0(T^kx),\eqsp x\in L.  
\end{equation} 
A linear functional $\sigma$ on $L$ is a summation if and only if the cocycle  
\eqref{eq:coc} is the coboundary of $\sigma$, see \eqref{eq:5}. 
Accordingly, $L$ admits summation 
if and only if the cocycle \eqref{eq:coc} is cohomologically trivial 
in the class of linear cochains.

Later on we assume that $A$ is provided with
a measure $\dmesns{\alpha}$, $\mes A=1$, and $f$ is a measure preserving 
transformation of $A$ into itself.  In this setting all cochains are
assumed measurable, and, accordingly, two cochains which coincide 
almost everywhere (a.e.) can be identified. (This is not necessary for our purposes.) 

The following lemma can be extracted from \cite{schmidt77} (see also \cite{moore80}, Section 5). 

\begin{lemma}
\label{lem:Sm}
Let $\psi(\alpha), \alpha\in A$, be a measurable function and let $\varepsilon > 0$.
Then there exist $M>0$ and a sequence of subsets $A_n\subset A$ such that
$\mes A_n >1-\varepsilon$ and 
\begin{equation}
\label{eq:bou}
|\psi(\alpha) - \psi(f^n\alpha)|\leq M,\eqsp \alpha\in A_n. 
\end{equation}
\end{lemma}
\begin{proof}  
  There is a subset $D$ such that $|\psi(\alpha)| \leq M/2, \alpha\in D$ and 
  $\mes(A\setminus D) < \varepsilon/2$. The inequality \eqref{eq:bou} is valid 
  on the intersection $A_n$ of $D$ with the preimage $f^{-n} D$.
  On the other hand, $\mes A_n> 1-\varepsilon$ since
\begin{equation*}  
  \mes(A\setminus A_n) \leq \mes(A\setminus f^{-n}D) + \mes(A\setminus D) =
  2\mes(A\setminus D) < \varepsilon.
\end{equation*} 
\end{proof} 

\begin{theorem}
\label{prop:5.6}
Let there exist a sequence of subsets
$B_n\subset A$, $\inf_n(\mes B_n) > 0$, and
\begin{equation}
\label{eq:5.26}
\inf_{\alpha\in B_n} \left|\sum_{k=0}^{n-1}\xi_0(f^k\alpha)\right|\rightarrow\infty,
\eqsp n\rightarrow\infty.
\end{equation}
Then c.e. \eqref{eq:5.8} has no measurable solutions.
\end{theorem}

\begin{proof}
We use Lemma \ref{lem:Sm} with 
$\varepsilon < \inf_n(\mes B_n)$ then $A_n\cap B_n\neq\emptyset$.
For $ \alpha\in A_n\cap B_n$ the equality 
\begin{equation*} 
\sum_{k=0}^{n-1} \xi_0(f^k\alpha) = \psi(\alpha) - \psi(f^n\alpha)
\end{equation*} 
yields
\begin{equation*} 
\inf_{\alpha\in B_n} \left|\sum_{k=0}^{n-1}\xi_0(f^k\alpha)\right| \leq M ,
\end{equation*} 
in contrary to  \eqref{eq:5.26}.
\end{proof}

Now we consider the space $L_1(A,\dmesns{\alpha})$ of Lebesgue 
integrable complex-valued functions. In this setting
the following Birkhoff-Khinchin 
ergodic theorem (see e.g. \cite{sinaietalergodictheory}, Ch.1) 
is our main tool.

\begin{theorem}
  \label{thm:BKT}
  Let $\phi\in L_1(A,\dmesns{\alpha})$. Then the limit
  \begin{equation*}
    \tilde{\phi}(\alpha) = \lim_{m\rightarrow\infty}
    \frac{1}{m}\sum_{n=0}^{m-1} \phi(f^n\alpha)
   \end{equation*}
  exists for $\alpha\in A_\phi$ where $A_\phi$ is an $f$-invariant
  subset of $A$, $\mes(A\setminus A_\phi)=0$. The limit function
  $\tilde{\phi}$ is $f$-invariant, it belongs to $L_1(A,\dmesns{\alpha})$,
  and 
  \begin{equation*}
    \int\tilde{\phi}\dmes{\alpha} = \int\phi\dmes{\alpha}.
  \end{equation*}
\end{theorem}

Recall that $f$ is said to be \defin{ergodic} if every 
$f$-invariant measurable function is constant a.e.. In this case
\begin{equation}
  \label{eq:5.11}
  \lim_{m\rightarrow\infty} \frac{1}{m} \sum_{n=0}^{m-1}
  \phi(f^n\alpha) = \int\phi\dmes{\alpha}, \eqsp \alpha\in A_\phi, 
\end{equation} 
where $A_\phi$ may be not the same as before, but has the same properties.  
Later on we deal with $A_\phi$ from \eqref{eq:5.11}.

\begin{theorem}
  \label{thm:5.4}
  Let $f$ be ergodic, and let the function $\xi_0\in L_1(A,\dmesns{\alpha})$ 
be such that
  \begin{equation} 
    \label{eq:5.18}
 \int \xi_0\dmes{\alpha} = 0.
  \end{equation}  
  Then the orbital series \eqref{eq:5.5} is summable on $A_{\xi_0}$, hence a.e.. 
\end{theorem}

\begin{proof}
  By Theorem \ref{thm:5} it suffices to show that
  $\pi_0\not\in L_X$. Suppose to the contrary. Then
  \begin{equation}
  \label{eq:5.60}
    \sum_k \lambda_k \xi_0(f^n\alpha_k) = 1, \eqsp n\geq 0,
  \end{equation}
  for a finite set
  $\{(\alpha_1,\lambda_1),(\alpha_2,\lambda_2),\ldots\}$ with $\alpha_k\in
  A_{\xi_0}$, $\lambda_k\in\C$. This contradicts \eqref{eq:5.11} with $\phi=\xi_0$. 
  Indeed, by \eqref{eq:5.18} the averaging (in the sense of \eqref{eq:5.11}) over $n$ in \eqref{eq:5.60} yields 0 on the left 
  hand side, instead of 1 on the right. 
\end{proof}

\begin{remark} 
\label{rem:5.2} 
Obviously, for any measure preserving $f$
the $L_1$-solvability of \eqref{eq:5.8} implies that
$\xi_0$ belongs to $L_1(A,\dmesns{\alpha})$ and satisfies \eqref{eq:5.18}.
Moreover, the latter is necessary for the existence of a measurable solution
to \eqref{eq:5.8}, see \cite{anosov73}. However, it is not sufficient. For
the irrational rotations of the circle and continuous $\xi_0$ this was shown
in \cite{anosov73} with the references to some dynamical constructions due to Neumann and
Kolmogorov. (For another construction see \cite{gordon75}.) In \cite{lyubich80}
the nonexistence of measurable solutions was established by means of the   
Banach theorem on closed graph. (See \cite{belitskii98} for a generalization.) 
Also note that the measurable solutions may be not Lebesgue integrable 
\cite{anosov73}, \cite{kornfeld76}.
\end{remark}

\begin{remark}
\label{rem:fur}
For a multiplicative version of c.e. the absence of measurable solutions
was proven in \cite{furstenberg61} assuming that the known function in the equation is
not homotopic to a constant. For this reason the problem for the additive equation
cannot be reduced to the result of  \cite{furstenberg61} by exponetiating.
\end{remark} 

In \cite{kolmogorov25} Kolmogorov claimed (without any proof or 
heuristics) that \defin{if the trigonometric series}
\begin{equation}
\label{eq:5.20}
  \sin t + \sin 3t + \cdots + \sin 3^n t + \cdots, \eqsp t\in\R, 
\end{equation}
\defin{is summable, then "one can construct an effective example of a Lebesgue
nonmeasureable function"}. Formally, the last sentence sounds as  
\defin{"the sum (in the sense of a summation) of the series \eqref{eq:5.20} is nonmeasurable"}. 
This property was proven by Zygmund (\cite{zygmund59}, Ch. 5, Problem 26) for the series 
\begin{equation}
\label{eq:5.200} 
  \cos t + \cos 2t + \cdots + \cos 2^n t + \cdots, \eqsp t\in\R.
\end{equation}
Our general theory allows us to 
prove Kolmogorov's conjecture in the form: \defin{the series \eqref{eq:5.20} 
is summable a.e., and its sum is nonmeasurable}. The same is true 
for the series \eqref{eq:5.200}. (It is interesting that \eqref{eq:5.200}
turns into the nonsummable series $\pi_0$ at $t=0$.) Moreover, we prove

\begin{theorem}
\label{thm:5.80}
Let $q$ be an integer, $q\geq 2$. Then
\begin{enumerate} 
\item  For any 
$2\pi$-periodic function $\theta\in L_1 (0,2\pi)$ 
with zero mean value the series
\begin{equation} 
\label{eq:5.81}  
\theta (t)+ \theta (qt) + \cdots  + \theta (q^n t) + \cdots
\end{equation} 
is summable a.e. to a function $\psi(t)$. 
     
\item $\psi(t)$ satisfies the c.e.
\begin{equation} 
\label{eq:5.810} 
\psi(t) - \psi(qt) = \theta(t). 
\end{equation} 
     
\item Let $\theta$ be a trigonometric polynomial, 
\begin{equation} 
\label{eq:theta}
\theta(t) = \sum_{i=1}^{m}(a_i\cos\nu_i t + b_i\sin\nu_i t),
\end{equation}  
and let all ratios $\nu_i/\nu_j$ $(i> j)$ be not the powers of $q$.
Then  all solutions to the equation \eqref{eq:5.810} are nonmeasurable. 
\end{enumerate}
\end{theorem}

In particular, in (3) $\theta(t)$ can be any trigonometric polynomal 
of degree $<q$.
 
\begin{proof}   
The transformation $f_q:t\rightarrow qt$ (mod $2\pi$) is ergodic. 
Hence, (1) follows from 
Theorem \ref{thm:5.4}. Then Proposition \ref{prop:5.1} implies (2). 
To prove (3) we use Theorem \ref{prop:5.6}. 

Consider the sequence of trigonometric polynomials 
\begin{equation}
\label{eq:teta}
\theta_n(t)= \sum_{k=0}^{n-1}\theta(q^k t),\eqsp n\geq 1.
\end{equation}
The Fourier spectrum $\Omega_n$ of $\theta_n(t)$ is the union of the pairwise disjoint sets 
$\{q^k\nu_i\}_{i=1}^{m}, 0\leq k \leq n-1$. Accordingly, the summands in \eqref{eq:teta} 
are pairwise orthogonal. Moreover, they have the same $L_2$-norm, say $\tau$.
Therefore, the $L_2$-norm of $\theta_n(t)$ is equal to $\tau\sqrt{n}$.
On the other hand, the sets $\Omega_n$ are uniformly lacunar: 
there is $\kappa>1$ independent of $n$ such that $\shift{\omega}\geq \kappa\omega$ for 
all $\shift{\omega},\omega\in \Omega_n$, $\shift{\omega}>\omega$. Indeed, let 
$\shift{\omega}=q^k\nu_i$ and $\omega=q^l\nu_j$. Then either $\shift{\omega}\geq 2\omega$, 
or $q^{k-l}< 2\max\{\nu_j/\nu_i: 1\leq i,j\leq m\}$. In the second case the set of all 
possible differences $k-l$ is finite since, in addition,   
$q^{k-l}>\min\{\nu_j/\nu_i: 1\leq i,j\leq m\}$. Hence, the latter inequality 
can be strengthened by inserting of a factor $\kappa>1$ into the right hand side. This yields 
$\shift{\omega}>\kappa\omega$. (Obviously, $\kappa<2$ if the second case 
is nonempty, otherwise, $\kappa=2$.)  

By virtue of the established properties of $\theta_n(t)$ 
there are some numbers $\gamma, \delta >0$ 
(depending on $\kappa$ only) such that the measure of every set
\begin{equation*} 
B _n = \left\{t: \left|\sum_{k=0}^{n-1}\theta_n(t)\right|\geq \gamma\sqrt{n}\right\} 
\end{equation*} 
is greater than $\delta$, see \cite{zygmund59} (Ch.5, Th. 8.25). 
Thus, Theorem \ref{prop:5.6} is applicable. 
\end{proof} 

\begin{cor}
\label{cor:trig}
If $\theta$ is a trigonometric polynomial such that the c.e. \eqref{eq:5.810} 
has a measurable solution $\psi$. Then $\psi(t)$ is a trigonometric polynomial a.e..
\end{cor}

\begin{proof}
By Theorem \ref{thm:5.80} there is $\nu_i\equiv 0$(mod $q$) in 
\eqref{eq:theta}. Let $\nu = \nu_l = q\mu$ be the maximum of such $\nu_i$. 
We will argue by induction on $\nu$. Consider  
\begin{equation*}
\tilde{\theta}(t) = \theta(t) + a_l(\cos\mu t - \cos\nu t)+ b_l(\sin\mu t - \sin\nu t). 
\end{equation*}
Accordingly, we introduce 
\begin{equation*}
\tilde{\psi}(t) = \psi(t) + a_l\cos\mu t + b_l\sin\mu t , 
\end{equation*}
so that $\tilde{\psi}(t) - \tilde{\psi}(qt) = \tilde{\theta}(t)$. If $\tilde{\theta} = 0$ 
then $\tilde{\psi}$  is a trigonometric polynomial a.e. since  
there is $\tilde{\nu}<\nu$ in the role of $\nu$ for $\tilde{\theta}$. 
If $\tilde{\theta} = 0$ then $\tilde{\psi}$ is a constant a.e. by ergodicity. As a result,  
$\psi$ is a trigonometric polynomial a.e. in any case.
\end{proof}

Now we can explicitly describe all the "trigonometric coboundaries" $\theta$.
\begin{theorem}
\label{thm:trcb}
A general form of the trigonometric coboundaries is 
\begin{equation}
\label{eq:genf} 
\theta(t) = \sum_{p\in I_d}\sum_{i=0}^{i_{p,d}}\left(a_{p,i}\cos pq^it + b_{p,i}\sin pq^it\right) 
\end{equation} 
\end{theorem}
\defin{where} $d\geq 1$, $I_d =\{p:1\leq p\leq d,\eqsp p\not\equiv 0$(mod $q$)\}, 
$i_{p,d} = \min\{i: pq^i > d\}$, 
\defin{and the coefficients satisfy}  
\begin{equation}
\label{eq:coef} 
\sum_{i=0}^{i_{p,d}}a_{p,i} =0,\eqsp \sum_{i=0}^{i_{p,d}}b_{p,i}=0.   
\end{equation}  
 
\begin{proof}
By substitution of  
\begin{equation*}
\psi(t)=\sum_{j=1}^d (h_j\cos jt + g_j\sin jt)
\end{equation*} 
into \eqref{eq:5.810} we obtain \eqref{eq:genf} with 
\begin{equation} 
\label{eq:dif}
a_{p,0}=h_p,\eqsp  a_{p,i}=h_{pq^i}-h_{pq^{i-1}}\eqsp (1\leq i\leq i_{p,d}-1), \eqsp  a_{p,i_{p,d}}=-h_{pq^{i_{p,d}-1}}
\end{equation} 
and similar formulas for $b_{p,i}$. The relations \eqref{eq:coef} follow from \eqref{eq:dif} by summation. 
This calculation is invertible since the representation $j=pq^i$ with $p\in I_{p,d}$ and $0\leq i\leq i_{p,d}-1$ 
is unique for every $j$, $1\leq j\leq d$. 
\end{proof}

In conclusion we return to Proposition \ref{prop:5.1} and inverse it as follows.
\begin{theorem}  
  \label{thm:5.1} 
  Let $f$ be ergodic, and let $\psi\in L_1(A,\dmes{\alpha})$ be a 
  solution of c.e. \eqref{eq:5.8} for $\alpha\in
  A^0$ where $A^0$ is an $f$-invariant subset of $A$, $\mes(A\setminus
  A^0)=0$. Then the formula 
  \begin{equation} 
    \label{eq:5.13}
    \sigma(X(\alpha)) = \psi(\alpha) - \int\psi\dmes{\alpha}, \eqsp
    \alpha\in A_1 =A^0\cap A_\psi . 
  \end{equation} 
determines a summation $\sigma$ of the orbital series \eqref{eq:5.5}.
\end{theorem}
  Let us emphasize that the set $A_1$ is $f$-invariant and
$\mes(A\setminus A_1) = 0$, so the series \eqref{eq:5.5} is summable a.e.. 
\begin{proof}    
  With any constant $c$ the function $\psi+c$ is also a solution of
  \eqref{eq:5.8} on $A^0$. In particular, such is 

\begin{equation*} 
    \hat{\psi}(\alpha) = \psi(\alpha) - \int\psi\dmes{\alpha},
  \end{equation*} 
  so \eqref{eq:5.13} can be rewritten as
  \begin{equation} 
    \label{eq:5.15} 
    \sigma(X(\alpha)) = \hat{\psi}(\alpha), \eqsp \alpha\in A_1,
  \end{equation}
  with 
  \begin{equation} 
    \label{eq:5.16} 
    \int\hat{\psi}\dmes{\alpha} = 0.
  \end{equation} 
  Formula \eqref{eq:5.15} correctly defines $\sigma (X(\alpha)), \alpha\in A_1$, if
\begin{equation*} 
  X(\alpha_1)= X(\alpha_2)\Rightarrow \hat{\psi}(\alpha_1)=\hat{\psi}(\alpha_2).
\end{equation*} 
Moreover, it can be extended linearly as long as 
  \begin{equation} 
    \label{eq:5.17} 
    \sum_k \lambda_k X(\alpha_k) = 0 \Rightarrow \sum_k
    \lambda_k\hat{\psi}(\alpha_k) = 0
  \end{equation} 
  for all finite sets $\{(\alpha_1,\lambda_1),  
  (\alpha_2,\lambda_2),\ldots\}$ with $\alpha_k\in A_1$,
  $\lambda_k\in\C$. 
  The resulting $\sigma$ is indeed a summation of $X(\alpha)$ on $A_1$ since
  \begin{equation*} 
    \sigma(X(\alpha)) - \sigma(T(X(\alpha))) = \hat{\psi}(\alpha) -
    \hat{\psi}(f\alpha) = \xi_0(\alpha).
  \end{equation*}
 It remains to prove the implication \eqref{eq:5.17}. 

The hypothesis in \eqref{eq:5.17} can be rewritten as
  \begin{equation*} 
    \sum_k \lambda_k\xi_0(f^l\alpha_k) = 0, \eqsp l\geq 0,
  \end{equation*}  
or, equivalently, as 
  \begin{equation*}
    \sum_k \lambda_k\hat{\psi}(f^l\alpha_k) - \sum_k \lambda_k
    \hat{\psi}(f^{l+1}\alpha_k) = 0, \eqsp l\geq 0. 
  \end{equation*} 
The sum of these equalities over $0\leq l\leq{n-1}$ yields 
  \begin{equation}
   \label{eq:v}
    \sum_k \lambda_k\hat{\psi}(\alpha_k) =  \sum_k \lambda_k
    \hat{\psi}(f^n\alpha_k), \eqsp n\geq 0. 
  \end{equation} 
 By \eqref{eq:5.16} the averaging over $n$ in \eqref{eq:v} 
 yields the conclusion in \eqref{eq:5.17}.
\end{proof} 

\begin{remark} 
\label{rem:fu}
Without any assumption on $\psi$ the c.e. \eqref{eq:5.8} is solvable if and only if 
$s_n(\alpha)=0$ for all $\alpha\in A, n\geq 1$, such that $f^n\alpha = \alpha$. The 
necessity of this condition is obvious. To the converse we introduce the equivalence 
relation on $A: f^m\beta = f^n\alpha$ for some $m,n$ (depending on $\alpha, \beta$).
It suffices to solve \eqref{eq:5.8} separately on each class of this equivalence, say, 
a class of an $\alpha$. To this end we determine  $\psi(f^n\alpha) =\psi(\alpha) - s_n(\alpha),n\geq 1$, 
and then  $\psi(\beta) = \psi(f^n\alpha) + s_m(\beta)$ as long as $f^m\beta = f^n\alpha$.    
It easy to show that $\psi$ is correctly defined and satisfies \eqref{eq:5.8}. For  
preperiodic $f$ an explicit solution has been given in \cite{belitskii98}.     
\end{remark}

\defin{Acknoledgement}. I am grateful to a referee for valuable remarks.

\bibliography{series}

\providecommand{\bysame}{\leavevmode\hbox to3em{\hrulefill}\thinspace}
\providecommand{\MR}{\relax\ifhmode\unskip\space\fi MR }
\providecommand{\MRhref}[2]{%
  \href{http://www.ams.org/mathscinet-getitem?mr=#1}{#2}
}
\providecommand{\href}[2]{#2}
\begin{thebibliography}{10}

\bibitem{anosov73}
D.~V. Anosov, \emph{On an additive functional homological equation connected
  with an ergodic rotation on the circle}, Math. USSR - Izv. \textbf{7} (1973),
  1257--1271.

\bibitem{belitskii98}
G.~R. Belitskii and Y.I. Lyubich, \emph{On the normal solvability of
  cohomological equations on compact topological spaces}, Oper. Th., Adv. and
  Appl. \textbf{103} (1998), 75--87.

\bibitem{furstenberg61}
H.~Furstenberg, \emph{Strict ergodicity and transformation of the torus},
  Amer.J.Math., \textbf{83} (1961), 573--601.

\bibitem{gordon75}
A.Ya. Gordon, \emph{Sufficient condition for unsolvability of the additive
  functional homological equation connected with the ergodic rotation of a
  circle}, Funct. Anal. Appl. \textbf{9} (1975), no.~4, 334--336.

\bibitem{gottschalk55}
W.~H. Gottschalk and G.~A. Hedlund, \emph{Topological dynamics}, AMS Coll.
  Publ., vol.~36, AMS, Providence, R. I., 1955.

\bibitem{hardy49}
G.~H. Hardy, \emph{Divergent series}, Univ. Press, Oxford, 1949.

\bibitem{hardy74}
\bysame, \emph{Collected papers}, vol.~6, Clarendon Press, Oxford, 1974.

\bibitem{katok95}
A.~Katok and B.~Hasselblatt, \emph{Introduction to the modern theory of
  dynamical systems}, Encycl. of Math. and Appl., vol.~54, Cambridge Univ.
  Press, Cambridge, 1995.

\bibitem{kirillov67}
A.~A. Kirillov, \emph{Dynamical systems, factors and group representations},
  Russian Math. Surveys \textbf{22} (1967), no.~5 (137), 63--75.

\bibitem{kolmogorov25}
A.~N. Kolmogorov, \emph{Sur la possibilit\'e de la d\'efinition g\'en\'erale de
  la d\'eriv\'ee, de l'int\'egrale et de la sommation de s\'eries divergentes},
  C. R. Acad. Sci. \textbf{180} (1925), 362--364.

\bibitem{kornfeld76}
I.~Kornfeld, \emph{On the additive homological equation}, Funct. Anal. Appl.
  \textbf{10} (1976), no.~2, 73--74.

\bibitem{sinaietalergodictheory}
I.~P. Kornfeld, S.~V. Fomin, and Ya.~G. Sina{\u\i}, \emph{Ergodic theory},
  Fundamental Principles of Mathematical Sciences, vol. 245, Springer-Verlag,
  New York, 1982.

\bibitem{lyubich80}
Yu.~I. Lyubich, \emph{Method of closed graph for the additive homological
  equation on the circle}, Yaroslavl Univ., Yaroslavl, 1980, in Russian.

\bibitem{lyubich92}
Yu.I. Lyubich, \emph{Linear functional analysis}, Functional Analysis (N.~K.
  Nikolskij, ed.), Encycl. of Math. Sci., vol.~19, Springer-Verlag, Berlin,
  1992.

\bibitem{moore80}
C.~C. Moore and K.~Schmidt, \emph{Coboundaries and homomorphisms for
  nonsingular actions and a problem of {H}. {H}elson}, Proc. London Math. Soc.
  (3) \textbf{40} (1980), no.~3, 443--475.

\bibitem{schmidt77}
K.~Schmidt, \emph{Cocycles on ergodic transformation groups}, Macmillan
  Lectures in Math., vol.~1, Macmillan, Delhi, 1977.

\bibitem{zygmund59}
A.~Zygmund, \emph{Trigonometric series}, vol.~1, Cambridge Univ. Press, 1959.

\end{thebibliography}

  Address:

  {\it Department of Mathematics, 

    Technion, 32000, 

    Haifa, Israel}

  \smallskip
  email: {\it lyubich@tx.technion.ac.il}  

\end{document}